\documentclass[11 pt,leqno]{article}
\title{
On solvable subgroups of automorphism groups of right-angled Artin groups
}
\author{Matthew B. Day }
\date{October 25, 2009}
\usepackage{slashbox}
\usepackage{amsmath}
\usepackage{amsfonts}
\usepackage{amssymb}
\usepackage{amsthm}
\usepackage{latexsym}
\usepackage{epsfig}
\usepackage{xypic}
\usepackage{verbatim}
\usepackage{mathrsfs}

\theoremstyle{plain} \newtheorem{theorem}{Theorem}[section]
\theoremstyle{plain} \newtheorem{proposition}[theorem]{Proposition}
\theoremstyle{plain} \newtheorem{lemma}[theorem]{Lemma}
\theoremstyle{plain} 
\theoremstyle{plain} \newtheorem{corollary}[theorem]{Corollary}
\theoremstyle{plain} 
\theoremstyle{plain} 
\theoremstyle{plain} 
\theoremstyle{plain} 
\theoremstyle{plain} 
\theoremstyle{remark} 
\theoremstyle{definition} \newtheorem{remark}[theorem]{Remark}
\theoremstyle{definition} \newtheorem{definition}[theorem]{Definition}
\theoremstyle{definition} \newtheorem{theorem-definition}[theorem]{Theorem-Definition}
\theoremstyle{definition} 
\theoremstyle{plain} 
\theoremstyle{plain} 
\theoremstyle{plain} 
\theoremstyle{plain} 

\numberwithin{equation}{section}

\DeclareMathOperator{\Out}{Out}
\DeclareMathOperator{\Aut}{Aut}
\DeclareMathOperator{\length}{length}
\DeclareMathOperator{\SL}{SL}

\newcommand\R{\mathbb{R}}

\newcommand\Z{\mathbb{Z}}

\newcommand\into\hookrightarrow

\DeclareMathOperator{\depth}{depth}
\DeclareMathOperator{\st}{st}
\DeclareMathOperator{\lk}{lk}
\newcommand\AAG{\Aut A_\Gamma}

\providecommand{\abs}[1]{\lvert#1\rvert}
\newcommand\OAG{\Out A_\Gamma}
\DeclareMathOperator{\GL}{GL}

\begin{document}
\maketitle

\begin{abstract}
For any right-angled Artin group, we show that its outer automorphism group contains either a finite-index nilpotent subgroup or a nonabelian free subgroup.
This is a weak Tits alternative theorem.
We find a criterion on the defining graph that determines which case holds.
We also consider some examples of solvable subgroups, including one that is not virtually nilpotent and is embedded in a non-obvious way.
\end{abstract}

\section{Introduction and Background}
\subsection{Introduction}
Let $A_\Gamma$ be the right-angled Artin group of a finite simplicial graph $\Gamma$ with vertex set $X$, \emph{i.e.} the group with presentation
\[A_\Gamma=\langle X| \{xy=yx | \text{ $x$ is adjacent to $y$ in $\Gamma$}\}\rangle.\]
In this note we find a combinatorial condition on the graph $\Gamma$ that indicates whether the outer automorphism group $\OAG$ of $A_\Gamma$ contains a non-abelian free group.
This extends a result of Gutierrez--Piggott--Ruane~\cite[Theorem~1.10]{gpr} which gives a condition for a particular subgroup of $\OAG$ to be abelian.
In fact, our theorem indicates a dichotomy: either $\Out A_\Gamma$ contains a nonabelian free subgroup or $\Out A_\Gamma$ is virtually nilpotent.
This is a weak Tits alternative theorem.
A true Tits alternative theorem would consider all subgroups of $\OAG$; Charney--Vogtmann~\cite{cvsubquot} recently proved such a theorem for a large class of right-angled Artin groups.

Automorphism groups of right-angled Artin groups are sometimes described as intermediate between automorphism groups of free groups and integer general linear groups, which are extreme examples.
This result is a first attempt to discern the cases where this idea seems reasonable, since $\Aut F_n$ and $\GL(n,\Z)$ both have nonabelian free subgroups for every $n>1$.

To state the theorem, we review some notions on graphs.
Recall that the \emph{link} $\lk(x)$ of a vertex $x\in\Gamma$ is the set of vertices adjacent to $x$, and the \emph{star} $\st(x)$ is $\lk(x)\cup\{x\}$.
\emph{Domination} is a useful relation that was considered by Servatius~\cite[Section~IV]{servatius}:
\begin{definition}
For $x,y\in\Gamma$, say $y$ \emph{dominates} $x$ if $\lk(x)\subset\st(y)$; denote this by $y\geq x$.
Say $x$ and $y$ are \emph{domination equivalent} if $x\leq y$ and $y\leq x$; denote this by $x\sim y$.
\end{definition}
Next we consider the notion of a \emph{separating intersection of links} defined by Gutierrez--Piggott--Ruane~\cite[Definition~1.9]{gpr}.
\begin{definition}
The graph $\Gamma$ has a \emph{separating intersection of links} if there are two vertices $x,y\in\Gamma$ such that \emph{(1)} $x$ is not adjacent to $y$ and \emph{(2)} there is a connected component of $\Gamma \backslash (\lk(x)\cap\lk(y))$ not containing $x$ or $y$.
\end{definition}

Now we state our main result.
\begin{theorem}\label{th:theorem}
Consider the following conditions on a graph $\Gamma$:
\begin{enumerate}
\item \label{it:sim}
$\Gamma$ contains a domination-equivalent pair of vertices.
\item \label{it:sil}
$\Gamma$ contains a separating intersection of links.
\end{enumerate}
If either condition holds, then $\OAG$ contains a nonabelian free subgroup.
If both conditions fail, then $\OAG$ is virtually nilpotent.
\end{theorem}

\begin{proof}[Proof of Theorem~\ref{th:theorem}]
If either condition holds, then $\OAG$ contains a nonabelian free subgroup, by Lemmas~\ref{le:equivfree} and~\ref{le:silfree} below.
If both conditions fail, then Proposition~\ref{pr:vnilclass} below produces a finite-index nilpotent subgroup.
\end{proof}

In Definition~\ref{de:depth}, we define a number $\depth(\Gamma)$ that can be read off of the graph $\Gamma$.
If $\OAG$ is virtually nilpotent, then a certain natural finite-index subgroup of $\OAG$ turns out to be nilpotent of class $\depth(\Gamma)$. 
Further, every finite-index nilpotent subgroup of $\OAG$ has nilpotence class at least $\depth(\Gamma)$.
See Proposition~\ref{pr:vnilclass} below for details.

In Section~\ref{ss:examples}, we consider a few other conditions that imply $\OAG$ contains a nonabelian free subgroup.
Then we construct examples of graphs $\Gamma$ with $\OAG$ containing finite-index nilpotent subgroups of arbitrary nilpotence class.

The corollary below follows from Theorem~\ref{th:theorem} by standard arguments.
\begin{corollary}\label{co:solvimpliesvnilp}
If $\Gamma$ is a graph such that $\OAG$ has a solvable, finite-index subgroup, then every solvable subgroup of $\OAG$ is virtually nilpotent.
\end{corollary}

Of course, for $n\geq 3$ the group $\GL(n,\Z)$ contains examples of solvable subgroups of infinite index that are not virtually solvable.
Given Corollary~\ref{co:solvimpliesvnilp}, one might conjecture that a solvable subgroup of $\OAG$ that is not virtually nilpotent must be essentially contained in an embedded copy of $\GL(n,\Z)$.
In Section~\ref{se:extra}, we produce an example where this is not the case.

\subsection{Background}
We will use the following four classes of automorphisms.
The \emph{inversion} of $x\in X$ is the automorphism sending $x$ to $x^{-1}$ and fixing $X-\{x\}$.
If $\pi$ is an automorphism (a symmetry) of the graph $\Gamma$, then the \emph{graphic automorphism} of $\pi$ is the automorphism sending $x$ to $\pi(x)$ for each $x\in X$.
If $x\in X\cup X^{-1}$ and $Y$ is a connected component of $\Gamma-\st(x)$, the \emph{partial conjugation} of $Y$ by $x$ is the automorphism sending $y$ to $x^{-1}yx$ for each $y\in Y$ and fixing $X-Y$. 
Denote this automorphism by $c_{x,Y}$.
If $x\in X\cup X^{-1}$ and $y\in X$ are distinct and $x\geq y$, then the \emph{transvection} of $y$ by $x$ is the automorphism sending $y$ to $yx$ and fixing $X-\{y\}$.
Denote this automorphism by $\tau_{x,y}$.
Sometimes we will refer to the automorphism just defined as the \emph{right transvection}, and refer to its conjugate by the inversion in $y$ as the \emph{left transvection}.
The \emph{multiplier} of a transvection $\tau_{x,y}$ is $x$ and the multiplier of a partial conjugation $c_{y,Y}$ is $y$.
Servatius defined these automorphisms and showed that they are well defined in~\cite[Section~IV]{servatius}.
Laurence~\cite{laurence} proved the following, which was a conjecture of Servatius.
\begin{theorem}[Laurence~\cite{laurence}]\label{th:laurence}
The finite set of all transvections, partial conjugations, inversions, and graphic automorphisms is a generating set of $\AAG$.
\end{theorem}
Of course the images of these generators form a finite generating set for $\OAG$.
Since we are working in $\OAG$, in this paper we will demand that partial conjugations are \emph{not} inner automorphisms.
Specifically, whenever we declare that $C_{y,Y}$ is a partial conjugation with multiplier $y$, we also assume $Y$ and $\Gamma\backslash (\st(y)\cup Y)$ are both nonempty.

\subsection{Acknowledgments}
I would like to thank the organizers of the 2009 International Conference on Geometric and Combinatorial Methods in Group Theory and Semigroup Theory, where I started considering the problem in this paper.
I am grateful to Ruth Charney for commenting on an earlier version of this paper.
This research was done under the support of an N.S.F. Mathematical Sciences Postdoctoral Research Fellowship.

\section{Proof of the dichotomy}
\subsection{Conditions for free subgroups}
\begin{lemma}\label{le:equivfree}
If $\Gamma$ contains distinct $x,y$ with $x\sim y$, then $\OAG$ contains a nonabelian free subgroup.
\end{lemma}

\begin{proof}
Let $x,y\in\Gamma$ be vertices with $x\sim y$.
Let $G$ be the subgroup of $\OAG$ generated by the images of $\tau_{x,y}^2$ and $\tau_{y,x}^2$.
The vector space $H_1(A_\Gamma;\R)$ has a basis given by the vertex-set of $\Gamma$; note that $G$ leaves the $2$--dimensional subspace $V=\langle[x],[y]\rangle$ invariant.
Let $A\subset V$ be the set of $a[x]+b[y]$ with $\abs{a}>\abs{b}$ and let $B\subset V$ be the set of such vectors with $\abs{b}>\abs{a}$.
It is easy to see that $A$ and $B$ are nonempty, $(\tau_{x,y}^2)_*(B)\subset A$ and $(\tau_{y,x}^2)_*(A)\subset B$.
Then by the well-known Table-Tennis Lemma (see de La Harpe~\cite[II.B.24]{delaharpe}), we see that $G$ is free of rank $2$.
\end{proof}

\begin{lemma}
\label{le:silfree}
If $\Gamma$ contains a separating intersection of links, then $\OAG$ contains a nonabelian free subgroup.
\end{lemma}
\begin{proof}
Let $x,y,z\in\Gamma$ with $y$ not adjacent to $z$ and $\lk(y)\cap\lk(z)$ separating $x$ from both $y$ and $z$.
Let $Y$ be the component of $x$ in $\Gamma\backslash \lk(y)$ and let $Z$ be the component of $x$ in $\Gamma\backslash \lk{z}$.
The hypotheses imply that $y\notin Z$ and $z\notin Y$.
Let $\widetilde G$ be the subgroup of $\AAG$ generated by $c_{y,Y}$ and $c_{z,Z}$ and let $G<\OAG$ be its image.
Then $\widetilde G$ fixes $y$ and $z$, and therefore contains no nontrivial inner automorphisms.
Therefore the projection $\widetilde G\to G$ is an isomorphism.

Map $\widetilde G$ to the free group $F_2=\langle y,z\rangle$ by sending $\alpha\in \widetilde G$ to the unique $w\in \langle y,z\rangle$ with $\alpha(x)=w^{-1}xw$.
It is easy to see that this map is a homomorphism with $c_{y,Y}$ mapping to $y$ and $c_{z,Z}$ mapping to $z$.
Since an inverse homomorphism $F_2\to\widetilde G$ is easy to construct, we see that $\widetilde G$ and $G$ are free of rank $2$.
\end{proof}
%**********************************************************
\subsection{Conditions for virtual nilpotence}

\begin{lemma}\label{le:othercondition}
Suppose $x$, $y$ and $z$ are in $\Gamma$ such that $x$ is not adjacent to $y$, $\lk(x)$ separates $y$ from $z$ and $\lk(y)$ separates $x$ from $z$.
Then $\lk(x)\cap\lk(y)$ separates $x$ and $y$ from $z$ and therefore $\Gamma$ contains a separating intersection of links.
\end{lemma}

\begin{proof}
Suppose $\lk(x)\cap\lk(y)$ does not separate both $x$ and $y$ from $z$.
Then there is a shortest path from $z$ to $x$ or $y$ through $\Gamma\backslash (\lk(x)\cap\lk(y)$.
Starting from $z$, the first time this path hits $\lk(x)\cup\lk(y)$ must also be the last, or else there would be a shorter path.
Then the hypotheses imply that the point on the path in $\lk(x)\cup \lk(y)$ must also be in $\lk(x)\cap\lk(y)$, a contradiction.
\end{proof}

\begin{lemma}\label{le:nonadjdom}
Suppose $x,y\in \Gamma$ with $x$ not adjacent to $y$, we have $x\geq y$ and  $\st(y)$ separates $\Gamma$.
Then $\Gamma$ contains a separating intersection of links.
\end{lemma}

\begin{proof}
Let $z$ be in a component of $\Gamma\backslash\st(y)$ not containing $x$.
This means $\lk(y)$ separates $z$ from $x$.
Since $x\geq y$, we know $\lk(x)\cap\lk(y)=\lk(y)$, so $\lk(x)\cap\lk(y)$ separates $z$ from $x$ and $y$.
Therefore $\Gamma$ contains a separating intersection of links.
\end{proof}

\begin{lemma}\label{le:generalthing}
Suppose $\Gamma$ does not contain a separating intersection of links.
Suppose $\alpha$ and $\beta$ are automorphisms that are either partial conjugations or transvections (or one of each) and that $\alpha$ and $\beta$ fix each other's multipliers.
Then $\alpha$ and $\beta$ commute in $\AAG$.
\end{lemma}

\begin{proof}
Let $x$ be the multiplier of $\alpha$ and let $y$ be the multiplier of $\beta$.
If $x=y^{\pm1}$ or $x$ is adjacent to $y$ then $\alpha$ and $\beta$ commute, so assume $x$ and $y$ are distinct and not adjacent.
Suppose there is some $z\in\Gamma$ such that neither $\alpha$ nor $\beta$ fixes $z$.

Suppose $x$ is adjacent to $z$.
Then $\alpha$ is a transvection (partial conjugations fix the links of their multipliers) and $x\geq z$.
If $y\geq z$, then $y$ is adjacent to $x$, counter to our assumption.
If $y\not\geq z$, then $\beta$ is a partial conjugation and $y$ is not adjacent to $z$.
This implies that $x$ and $z$ are in the same connected component of $\Gamma\backslash \lk(y)$, meaning that $\beta$ cannot fix $x$ and change $z$.
This contradiction implies that $x$ is not adjacent to $z$, and similarly, that $y$ is not adjacent to $z$.

Suppose $\lk(x)$ does not separate $y$ from $z$.
Then $x\not\geq z$ and $\alpha$ is a partial conjugation.
However, in that case $\alpha$ cannot fix $y$ and change $z$.
So $\lk(x)$ separates $y$ from $z$, and similarly $\lk(y)$ separates $x$ from $z$.
Then by Lemma~\ref{le:othercondition}, we have that $\Gamma$ does contain a separating intersection of links, which is a contradiction.
From this we deduce that for each $z\in\Gamma$, either $\alpha$ fixes $z$ or $\beta$ fixes $z$.
This is enough to deduce that $\alpha$ and $\beta$ commute.
\end{proof}

In the case that $\alpha$ and $\beta$ are partial conjugations, the following lemma is a special case of Theorem~1.10 from Gutierrez--Piggott--Ruane~\cite{gpr}.
\begin{lemma}\label{le:commute}
Suppose $\Gamma$ does not contain a separating intersection of links.
Let $\alpha$ be a partial conjugation.
Suppose that $\beta$ is a transvection fixing the multiplier of $\alpha$, or that $\beta$ is a partial conjugation (not necessarily fixing the multiplier of $\alpha$).
Then the images of $\alpha$ and $\beta$ commute in $\OAG$.
\end{lemma}
\begin{proof}
Suppose the multipliers of $\alpha$ and $\beta$ are distinct (otherwise $\alpha$ and $\beta$ commute).
Then possibly by multiplying $\alpha$ and $\beta$ by inner automorphisms, we may assume that $\alpha$ fixes the multiplier of $\beta$, and if $\beta$ is a partial conjugation, we may assume that $\beta$ fixes the multiplier of $\alpha$.
The lemma then follows from Lemma~\ref{le:generalthing}.
\end{proof}

\begin{lemma}\label{le:steinberg}
Suppose conditions~\eqref{it:sim} and~\eqref{it:sil} from Theorem~\ref{th:theorem} both fail.
Suppose $\alpha=\tau_{x,y}$ is a transvection and $\beta$ is a transvection or partial conjugation with multiplier $y^{\pm1}$.
Let $\gamma$ be any commutator of $\alpha$ or $\alpha^{-1}$ with $\beta$ or $\beta^{-1}$.
Then $\gamma$ is a transvection or partial conjugation with multiplier $x^{\pm1}$ in $\OAG$.
Further, if $\beta$ is a transvection acting on $z^{\pm1}$, then so is $\gamma$, and if $\beta$ is a partial conjugation acting on $Y\subset X$, then so is $\gamma$.
\end{lemma}

\begin{proof}
If $\beta$ is a transvection and doesn't fix $x$, then $x\sim y$ and condition~\eqref{it:sim} holds.
If $\beta$ is a partial conjugation, then up to an inner automorphism we may assume that it fixes $x$.
So assume $\beta$ fixes $x$.

We claim that $x$ is adjacent to $y$.
Let $z\in \Gamma$ be an element not fixed by $\beta$.
If $y$ is adjacent to $z$, then $\beta$ is a transvection and $y\geq z$.
Since $x\geq y$, this implies that $x\geq z$, and that $x$ is adjacent to $y$.
So suppose $y$ is not adjacent to $z$.
If $x$ is not adjacent to $y$, then $\lk(y)$ separates $x$ from $z$ since $\beta$ fixes $x$ but not $z$.
Then by Lemma~\ref{le:nonadjdom}, $\Gamma$ contains a separating intersection of links, contradicting the failure of condition~\eqref{it:sil}.
So $x$ is adjacent to $y$.
Then the lemma follows by a computation.
\end{proof}

\begin{definition}
A \emph{domination chain} in $\Gamma$ is a sequence of distinct vertices $x_1,\dotsc,x_m$ of $\Gamma$ such that $x_m\geq x_{m-1}\geq\dotsb\geq x_1$.
The \emph{length} of the domination chain $x_1,\dotsc,x_m$ is $m-1$.
The \emph{domination depth} of $x$ is the length of the longest domination chain with $x$ as the dominant member.
\end{definition}

\begin{definition}
A domination chain $x_m\geq \dotsb\geq x_1$ is \emph{star-separation preserving} if $\Gamma\backslash\st(x_1)$ has two components $Y_1$ and $Y_2$ such that $Y_i\not\subset\st(x_m)$ for $i=1,2$.
The \emph{star-separation depth} of $x\in\Gamma$ is 
\[1+\max_{x=x_m\geq\dotsb\geq x_1}\length(x_m\geq\dotsb\geq x_1)\]
where the maximum is taken over all star-separation-preserving domination chains.
\end{definition}

\begin{definition}\label{de:depth}
The \emph{depth} $\depth(x)$ of a vertex $x\in \Gamma$ is maximum of the domination depth of $x$ and the star-separation depth of $x$.
The \emph{depth} $\depth(\Gamma)$ of $\Gamma$ is the maximum depth of its vertices.
\end{definition}

\begin{proposition}\label{pr:vnilclass}
The subgroup $N$ of $\OAG$ generated by transvections and partial conjugations is finite index in $\OAG$.
If the conditions from Theorem~\ref{th:theorem} both fail, then $N$ is nilpotent of class $\depth(\Gamma)$.
Further, every finite-index nilpotent subgroup of $\OAG$ has nilpotence class at least $\depth(\Gamma)$.
\end{proposition}

\begin{proof}
Let $S$ be the finite subset of $\OAG$ consisting of the identity, the images of transvections (both right and left) and partial conjugations, and their inverses.
Let $N$ be the subgroup generated by $S$ and let $P$ be the finite subgroup generated by images of inversions and graphic automorphisms in $\OAG$.
Note that $P$ normalizes $N$ (since conjugation by $P$ leaves $S$ invariant).
By Laurence's theorem (Theorem~\ref{th:laurence}), $\OAG=P N$ and therefore $N\lhd \OAG$.
By a classical group isomorphism theorem, $\OAG/N\cong P/(P\cap N)$ and therefore $N$ is finite-index in $\OAG$.
(In fact, the failure of condition~\eqref{it:sim} implies that $N\cap P=1$ and $\OAG\cong P\ltimes N$, as can be seen from the presentation for $\AAG$ in Day~\cite[Theorem~2.7]{day}.)

Let $k=\depth(\Gamma)$.
Let $S_0=\{1\}$, and let $S_i$ be the union of $\{1\}$ with the set of tranvsections $\tau_{x,y}$ with $\depth(x)-\depth(y)\geq k-i+1$ (and left transvections satisfying the same condition) and partial conjugations $c_{y,Y}$ with $\depth(y)\geq k-i+1$ for $i=1,\dotsc,k$.
The $S_i$ are nested and $S$ is $S_k$.

Let $\alpha\in S_i$ and $\beta\in S_j$, for $1\leq i,j\leq k$.
By Lemmas~\ref{le:generalthing},~\ref{le:commute} and~\ref{le:steinberg}, we see that if $[\alpha,\beta]$ is nontrivial, then $i+j>k+1$ and $[\alpha,\beta]$ is a member of $S_{i+j-k-1}$.
Since $i,j\leq k$, we have that $i+j-k-1<i,j$.
This is enough to deduce that $N$ is nilpotent of class at most $k$.

Select $x_k\in\Gamma$ with $\depth(x_k)=k$.
By definition, there is a domination chain $x_k\geq\dotsb\geq x_1$ in $\Gamma$, such that either $x_1$ dominates a vertex $x_0$, or $\Gamma\backslash\st(x_1)$ has two components $Y_1,Y_2$ with $Y_i\not\subset\st(x_k)$ for $i=1,2$ (depending on whether $\depth(x_k)$ is the domination depth or the star-separation depth, respectively).
In the first of these cases, let $\alpha_1$ denote the transvection $\tau_{x_1,x_0}$, and in the second of these cases, let $\alpha_1$ denote the partial conjugation $c_{x_1,Y_1}$.
For $i=2,\ldots, k$, let $\alpha_i$ be the transvection $\tau_{x_i,x_{i-1}}$.
Then by Lemma~\ref{le:steinberg}, the element 
\[[\dotsb[[\alpha_1,\alpha_2],\alpha_3],\dotsc,\alpha_{k}]\in\OAG\]
 is either a transvection $\tau_{y,x_0}$ or a partial conjugation $c_{y,Y}$, where $y=x_k^{\pm1}$ and $Y=Y_1\backslash\st(y)$.
If it is $\tau_{y,x_0}$, it is obviously nontrivial in $\OAG$.
If it is $c_{y,Y}$, it is nontrivial in $\OAG$ since there is an element of $Y_1\backslash \st(y)$ that is conjugated and an element of $Y_2\backslash \st(y)$ that is not conjugated.
So the nilpotence class of $N$ equals~$\depth(\Gamma)$.

Now suppose that $N''$ is a nilpotent, finite-index subgroup of $\OAG$.
Then $N''$ intersects $N$ in a finite index subgroup $N'$.
Each of the $\alpha_1,\dotsc,\alpha_k$ from the previous paragraph is of infinite order.
Since $N'$ is finite index in $N$, the intersection $N'\cap\langle \alpha_i\rangle$ is finite index in $\langle\alpha_i\rangle$ for each $i$.
In particular, each $N'\cap\langle \alpha_i\rangle$ is nontrivial.
So we have $a_1,\dotsc,a_k\in\Z$ with $\alpha_i^{a_i}\in N'$ for each $i$.
Then by the same reasoning as in the previous paragraph, we see that
\[[\dotsb[[\alpha_1^{a_1},\alpha_2^{a_2}],\alpha_3^{a_3}],\dotsc,\alpha_{k}^{a_k}].\]
 is nontrivial.
From this, we see the nilpotence class of $N'$ is also $\depth(\Gamma)$, and the nilpotence class of $N''$ is at least $\depth(\Gamma)$.
\end{proof}

The following needs no further proof.
\begin{corollary}
The group $\OAG$ is virtually abelian if and only if both conditions from Theorem~\ref{th:theorem} fail and $\depth(\Gamma)\leq 1$.
\end{corollary}

\begin{remark}
It has long been known that $\OAG$ is finite if and only if $\Gamma$ contains no pair of vertices $x, y$ with $x\geq y$ and $\Gamma$ contains no vertex $x$ with $\Gamma-\st(x)$ disconnected.
This is an easy corollary of Theorem~\ref{th:laurence}.
\end{remark}

\subsection{Examples}\label{ss:examples}

\begin{corollary}
The group $\OAG$ has a nonabelian free subgroup if any of the following conditions on $\Gamma$ hold:
\begin{itemize}
\item $\Gamma$ is disconnected.
\item $\Gamma$ contains a cut-vertex that breaks $\Gamma$ into three or more components.
\item $\Gamma$ contains non-adjacent vertices $x$ and $y$ with $x\geq y$ and $\st(y)$ separating $\Gamma$.
\item $\Gamma$ contains pairwise non-adjacent vertices $x$, $y$ and $z$ with $x\geq y\geq z$.
\end{itemize}
\end{corollary}

\begin{proof}
In each case we find a domination-equivalent pair of vertices or a separating intersection of links in $\Gamma$, and Theorem~\ref{th:theorem} implies the corollary.
The final condition is a special case of the second to last condition, which implies $\Gamma$ has a separating intersection of links by Lemma~\ref{le:nonadjdom}.

Now suppose that $\Gamma$ is disconnected.
If $\Gamma$ is edgeless, then any two vertices are domination equivalent.
Otherwise some component of $\Gamma$ has at least two vertices.
If each component of $\Gamma$ is a complete graph, then any two vertices in the same component are domination equivalent.
So we have some component of $\Gamma$ that contains two nonadjacent vertices.
Then $\Gamma$ contains a separating intersection of links (for $x$ and $y$ not adjacent, $\lk(x)\cap\lk(y)$ separates $x$ and $y$ from any vertex in another component).

Now suppose $\Gamma$ contains a cut-vertex $z$ that breaks $\Gamma$ into at least three components.
Without loss of generality we assume $\Gamma$ is connected.
If the valence of $z$ is less than $2$, then $\Gamma\backslash\{z\}$ has only one component.
If for for each pair of distinct $x,y\in\lk(z)$, either $x$ is adjacent to $y$ or $\lk(x)\cap\lk(y)$ contains two or more elements, then $\Gamma\backslash\{z\}$ has only one component.
Therefore $\Gamma$ contains distinct, non-adjacent vertices $x$ and $y$ with $\lk(x)\cap\lk(y)=\{z\}$.
Then $\Gamma\backslash(\lk(x)\cap\lk(y))$ has at least three components and $\Gamma$ has a separating intersection of links.
\end{proof}

\begin{proposition}
For each $k\geq 0$, there is a graph $\Gamma_k$ such that $\Out A_{\Gamma_k}$ contains a finite-index subgroup of nilpotence class $k$.
\end{proposition}

\begin{proof}
For each $k$, we will construct a graph $\Gamma_k$ with $\depth(\Gamma_k)=k$ and such that $\Gamma_k$ satisfies the hypotheses of Proposition~\ref{pr:vnilclass}.
We can take $\Gamma_0$ to be the graph with one vertex.

Now fix $k>0$.
For the vertex set of $\Gamma_k$, we will take a set of $2k+2$ vertices labeled as $x_0,\dotsc,x_k$, $y_0,\dotsc,y_k$.
Take the induced subgraph on $\{x_i\}_i$ to be the complete graph on $k$ vertices, and similarly for $\{y_i\}_i$.
Further, connect $x_i$ to $y_j$ by an edge if $i+j>k$.
These are the only edges of $\Gamma_k$.

Then $x_k\geq y_0$, $y_k\geq x_0$, and for $0\leq i<j\leq k$ we have $x_j\geq x_i$ and $y_j\geq y_i$.
Since $k>0$, these are the only pairs which satisfy the domination relation.
In particular, there are no domination-equivalent pairs.
There are no vertices whose stars separate $\Gamma$, so the star-separation depth of all vertices is trivial.
We compute all depths as equal to domination depths, and find $\depth(x_i)=\depth(y_i)=i$ for $i=0,\dotsc,k$.
Therefore $\depth(\Gamma_k)=k$.

The only non-adjacent pairs of vertices are $(x_i,y_j)$ and $(x_j,y_i)$ for $i+j\leq k$.
For such $i,j$, every element of $\Gamma_k$ is adjacent to either $x_i$ or $y_j$.
In particular, every element of $\Gamma_k\backslash(\lk(x_i)\cap\lk(y_j))$ has a path of length one to either $x_i$ or $y_j$ (and similarly for $x_j$ and $y_i$).
Therefore $\Gamma_k$ does not contain a separating intersection of links.
\end{proof}

\section{A non-nilpotent solvable subgroup}\label{se:extra}
Whenever $Y\subset X$ is a clique with $x\sim y$ for all $x,y\in Y$, the transvections of elements of $Y$ acting on each other generate an embedded copy of $\SL(|Y|,\Z)$ inside $\OAG$.
When we have such a copy of $\SL(n,\Z)$, say it is \emph{canonically embedded}.
Of course one can find non-virtually-nilpotent solvable subgroups of $\OAG$ inside canonically embedded copies of $\SL(n,\Z)$ for $n\geq 3$.
Given Corollary~\ref{co:solvimpliesvnilp}, one might conjecture that when $G<\OAG$ is solvable but not virtually nilpotent, there is $H<\OAG$ a canonically embedded copy of $\SL(n,\Z)$, such that $H\cap G$ is not virtually nilpotent.
However the following example is not of this type.
\begin{proposition}
Let $\Gamma$ be the graph on three vertices $\{a,b,c\}$ with a single edge from $a$ to $b$.
Let $G$ be the subgroup of $\OAG$ generated by the images of the elements $\{\tau_{a,c},\tau_{b,c},\tau_{a,b}\tau_{b,a}\}$.
Then $G$ is a solvable group and is not virtually nilpotent.

The intersection of $G$ with the unique canonically embedded copy of $\SL(n,\Z)$ in $\OAG$ is not virtually nilpotent.
\end{proposition}

\begin{proof}
It is apparent that $G$ does not contain any inner automorphisms, so $G$ is isomorphic to the subgroup of $\AAG$ generated by these generators.
Let $\alpha=\tau_{a,c}$, $\beta=\tau_{b,c}$, and let $\gamma=\tau_{a,b}\tau_{b,a}$.
Since $a$ commutes with $b$, we know that $\alpha$ commutes with $\beta$.
A computation shows that $\gamma\alpha\gamma^{-1}=\alpha^{2}\beta$ and $\gamma\beta\gamma^{-1}=\alpha\beta$.
It is easy to see that $\langle \alpha,\beta\rangle\cap\langle\gamma\rangle=1$.
From this we can see that $G$ is the semidirect product $\Z\ltimes \Z^2$, where $\Z$ acts on $\Z^2$ by the matrix $\binom{2\,\,1}{1\,\,1}$.
So $G$ is solvable.

On the other hand, for all $k>0$, the centralizer of $\gamma^k$ in $G$ is $\langle \gamma\rangle$.
Let $H$ be a finite index subgroup of $G$.
Then $H$ contains a positive power of $\gamma$ and an element of $G$ outside of $\langle \gamma\rangle$.
So $H$ has trivial center and is therefore not nilpotent.
In fact $G$ is isomorphic to a lattice in the $3$-dimensional Lie group $\mathrm{sol}$; see Thurston~\cite[Example~3.8.9]{thurston} for explanation.

The only canonically embedded copy of any $\SL(n,\Z)$ in $\OAG$ is generated by $\tau_{a,b}$ and $\tau_{b,a}$.
However, the intersection of $G$ with this subgroup is a copy of $\Z$.
\end{proof}

\bibliographystyle{amsplain}
\bibliography{nilraag}

\noindent
Department of Mathematics 253-37\\
California Institute of Technology\\
Pasadena, CA 91125\\
E-mail: {\tt mattday@caltech.edu}
\medskip

\end{document}